\documentclass[a4paper,11pt,twoside,reqno]{amsart}

\usepackage[utf8]{inputenc}
\usepackage[plainpages=false,pdfpagelabels=true]{hyperref}
\usepackage{amssymb,amsthm}
\usepackage[margin=1in]{geometry}

\newtheorem{Theorem}{Theorem}[section]
\newtheorem{Prop}[Theorem]{Proposition}
\newtheorem{Lem}[Theorem]{Lemma}

\newtheorem{Cor}[Theorem]{Corollary}
\theoremstyle{definition}
\newtheorem{Dfn}[Theorem]{Definition}
\newtheorem{Bem}[Theorem]{Remark}
\newtheorem{Bsp}[Theorem]{Example}
\newcommand{\tr}{\operatorname{Tr}}
\newcommand{\Ric}{\operatorname{Ric}}

\newcommand{\Scal}{\operatorname{Scal}}

\newcommand{\vol}{{\operatorname{vol}}}
\newcommand{\dv}{\text{ }dv}
\newcommand{\s}{\mathbb{S}}

\newcommand{\N}{\mathbb{N}}

\newcommand{\R}{\mathbb{R}}

\allowdisplaybreaks[1]

\renewcommand{\epsilon}{\varepsilon}

\numberwithin{equation}{section}

\usepackage{color}

\title{Remarks on constructing biharmonic and conformal-biharmonic maps to spheres}
\author{Volker Branding}
\address{University of Rostock, Institute of Mathematics\\
Ulmenstraße 69, 18057 Rostock, Germany}
\email{volker.branding@uni-rostock.de}
\thanks{The author gratefully acknowledges the support of the Austrian Science Fund (FWF) through the project "Geometric Analysis of Biwave Maps" (DOI: 10.55776/P34853)
}

\date{\today}
\subjclass[2010]{58E20}

\subjclass[2010]{58E20; 53C43}
\keywords{biharmonic map; conformal-biharmonic map; sphere}

\begin{document}

\begin{abstract}
Biharmonic and conformal-biharmonic maps are two fourth-order generalizations
of the well-studied notion of harmonic maps in Riemannian geometry.
In this article we consider maps into the Euclidean sphere and 
investigate a geometric algorithm that aims at rendering a given harmonic map either biharmonic or conformally biharmonic. 

For biharmonic maps we find that in the case of a closed domain the maximum principle imposes strong restrictions 
on our approach, whereas there is more flexibility when we have a non-compact domain and we highlight this difference by a number of examples.

Concerning conformal-biharmonic maps we show that our algorithm produces
explicit critical points for maps between spheres. Moreover, it turns out that
we do not get strong restrictions as we obtain for biharmonic maps, such that our algorithm might produce additional conformal-biharmonic maps between spheres beyond the ones found in this article. 
\end{abstract}

\maketitle

\section{Introduction and results}
Let us consider two Riemannian manifolds \((M,g),(N,h)\) and a map \(\phi\colon M\to N\) which we assume to be smooth throughout this manuscript. A particular interesting class of maps can be obtained by calculating the critical points of the energy of a map which is given by

\begin{align}
\label{eq:energy}
E(\phi)=\frac{1}{2}\int_M|d\phi|^2\dv.    
\end{align}

The critical points of \eqref{eq:energy} are characterized by the vanishing
of its tension field. i.e.
\begin{align}
\label{eq:harmonic-map}
0=\tau(\phi):=\tr\bar\nabla d\phi.
\end{align}
Here, \(\bar\nabla\) denotes the connection on the pull-back bundle \(\phi^\ast TN\).
The solutions of \eqref{eq:harmonic-map} are precisely \emph{harmonic maps}.
The harmonic map equation represents a second order semilinear elliptic partial differential equation and as of today many results concerning their existence and non-existence could be established. In this regard we would like to mention the seminal result of \emph{Eells \(\&\) Sampson} \cite{MR164306} which guarantees the existence of a harmonic map in each homotopy class of maps provided that \(M\) is closed and that \(N\) has non-positive curvature.
If the target manifold has positive curvature, i.e. in the case of a sphere, the question on the existence of harmonic maps is substantially more difficult. Here, we want to mention the seminal 
article of Smith \cite{MR391127} who introduced reduction techniques in order to turn the harmonic map equation into an ordinary differential equation. More details concerning this approach can be found in the book of Eells \(\&\) Ratto \cite{MR1242555}.

We would like to point out that harmonic maps are special if the domain is two-dimensional. In this case both the energy \eqref{eq:energy} as well as the harmonic map equation \eqref{eq:harmonic-map} are invariant under conformal transformations on the domain. Moreover, harmonic maps are energy critical from an analytic perspective in two dimensions.

For more details on the theory of harmonic maps we refer to the book
\cite{MR2044031}.

A fourth order generalization of harmonic maps that has received growing interest over the last years is the variational problem of \emph{biharmonic maps}. Here, the starting point is the \emph{bienergy functional} 
which is defined as follows
\begin{align}
\label{eq:bienergy}
E_2(\phi)=\frac{1}{2}\int_M|\tau(\phi)|^2\dv.
\end{align}
The critical points of \eqref{eq:bienergy} are characterized by the vanishing of the bitension field, i.e.
\begin{align}
\label{eq:bitension}
0=\tau_2(\phi):=\bar\Delta\tau(\phi)+\tr R^N(\tau(\phi),d\phi(\cdot))d\phi(\cdot),
\end{align}
where \(\bar\Delta\) represents the connection Laplacian on the pull-back bundle 
\(\phi^\ast TN\) and \(R^N\) the Riemannian curvature tensor on \(N\).
The critical points of the bienergy were first calculated by Jiang \cite{MR886529}. 
The biharmonic map equation \eqref{eq:bitension} comprises a semilinear equation of fourth order and the higher number of derivatives leads to significant problems in the mathematical analysis. In particular, powerful tools such as the maximum principle are no longer applicable in general.

A direct inspection of the biharmonic map equation \eqref{eq:bitension} shows that harmonic maps always provide a class of solutions. For this reason one is mostly interested in finding those biharmonic maps which are non-harmonic, these are usually called \emph{proper biharmonic}. In \cite{MR886529} Jiang also showed, with the help of the maximum principle, that in the case of \(M\) being closed and \(N\) having non-positive curvature every biharmonic map must be harmonic. 
Further such classification results can be found in \cite{MR4040175}.
Hence, one cannot find proper biharmonic maps in this geometric setup and, as a consequence, most attention is paid to constructing and classifying biharmonic maps to spheres. 
In \cite{MR2448058} it was conjectured that the only proper
biharmonic hypersurfaces of spheres are given by
\begin{enumerate}
    \item the small hypersphere \(\s^{n-1}(1/\sqrt{2})\),
    \item the generalized Clifford torus \(\s^{n_1}(1/\sqrt{2})\times\s^{n_2}(1/\sqrt{2})\), where
     \(n_1+n_2=n-1,n_1\neq n_2\).
\end{enumerate}
Here, \(\s^n\) represents the \(n\)-dimensional Euclidean sphere and
\(\s^n(r)\) denotes the Euclidean sphere of radius \(r\).

For the current status of research on biharmonic maps in Riemannian geometry
we refer to the recent book \cite{MR4265170}.

We want to mention, in sharp contrast to the energy of a map \eqref{eq:energy}, that the bienergy \eqref{eq:bienergy} is not invariant under conformal transformations of the domain metric in any dimension. Hence, from the point of view of conformal geometry, biharmonic maps might not be the correct fourth order generalization of harmonic maps.

Fortunately, there is a way to render the bienergy \eqref{eq:bienergy}
conformally invariant by augmenting it with two additional terms as follows
\begin{align}
\label{eq:bienergy-conformal}
E^c_2(\phi):=\frac{1}{2}\int_M \big(|\tau(\phi)|^2+\frac{2}{3}\Scal^M |d\phi|^2-2\tr\langle d\phi\left(\Ric^M(\cdot)\right), d\phi(\cdot)\rangle\big) \dv.
\end{align}
We call \eqref{eq:bienergy-conformal} the \emph{conformal bienergy functional}.
Here, \(\Scal^M\) denotes the scalar curvature of the domain while \(\Ric^M\) represents its Ricci curvature, considered as an endomorphism of \(TM\).

The critical points of the conformal bienergy functional \eqref{eq:bienergy-conformal} are those who satisfy
\begin{align}
\label{eq-c-biharmonic-intro}
0=\tau_2^c(\phi) & := \tau_2(\phi)-\frac{2}{3}\Scal^M\tau(\phi)+2\tr(\bar \nabla d\phi)(\Ric^M(\cdot),\cdot)+\frac{1}{3}d\phi(\nabla\Scal^M).
\end{align}
Solutions of \eqref{eq-c-biharmonic-intro} are called \emph{conformal-biharmonic maps}, or simply 
\emph{c-biharmonic maps}.

For more details and references on conformal-biharmonic maps we refer to the recent article \cite{BNO}.

As the equations for biharmonic and conformal-biharmonic maps are of fourth order
it is a challenging task to construct non-trivial solutions.
In order to overcome these technical difficulties it is thus reasonable
to start with a particular kind of map, 
of which we already have sufficient control, and then try to deform it to a biharmonic or conformal-biharmonic map.

The starting point of the calculations performed in this article are the following observations on biharmonic curves on spheres which were obtained in \cite{MR1863283},
see also \cite[Section 2.2]{MR4824947}.
Note that for a one-dimensional domain
the bienergy \eqref{eq:bienergy} coincides with the conformal bienergy \eqref{eq:bienergy-conformal} such that the equations for biharmonic and conformal-biharmonic maps are the same in this case.

There exist two classes of proper biharmonic curves parametrized by arc length on Euclidean spheres, which already exhaust all classes of proper biharmonic curves on this target
that have unit speed.
More precisely, we have:
\begin{enumerate}
    \item 
The proper biharmonic curve on \(\s^2\subset\R^3\) is given by the following parametrization 
\begin{align*}
\gamma(s)
=\big(\frac{1}{\sqrt{2}}\cos(\sqrt{2}s),\frac{1}{\sqrt{2}}\sin(\sqrt{2}s),
\frac{1}{\sqrt{2}}\big)
\end{align*}
which represents a curve of constant geodesic curvature \(k=1\) parametrized by arc length.

\item There exists a second class of proper biharmonic curves on \(\s^3\) given by
\begin{align}
\label{eq:biharmonic-curve-b}
\gamma(s)=\big(\frac{1}{\sqrt{2}}\sin(as),\frac{1}{\sqrt{2}}\cos(as),  
\frac{1}{\sqrt{2}}\sin(bs),\frac{1}{\sqrt{2}}\cos(bs)\big)
\end{align}
with \(a^2+b^2=2,a\neq b\). Note that these proper biharmonic curves have constant geodesic curvature \(k\) and torsion \(\tau\) and satisfy 
\(k^2+\tau^2=1\). 
The constraint \(a^2+b^2=2\) ensures that the curve is parametrized with respect to arc length while the condition \(a\neq b\) makes sure that the curve is not congruent to a geodesic.
\end{enumerate}
This classification of proper biharmonic curves on Euclidean spheres was first presented in \cite{MR1919374}.
For more details on biharmonic curves we refer to \cite{MR4542687}
and references therein.

Motivated by the first class of proper biharmonic curves on spheres
we will investigate when maps of the form
\begin{align}
\label{eq:definition-q-biharmonic}
q:=(\sin\alpha\cdot v,\cos\alpha),\qquad \alpha\in (0,\frac{\pi}{2}),  
\end{align}
where \(v\colon M\to\s^{n-1}\subset\R^{n}\) is a harmonic map 
can actually be biharmonic, i.e. a solution of \eqref{eq:bitension} or conformal-biharmonic, i.e. a solution of \eqref{eq-c-biharmonic-intro}.
More precisely, we want to understand if biharmonic maps from a higher-dimensional domain can be more flexible than in the case of curves. 

From a geometric point of view we know that for \(\alpha=\frac{\pi}{2}\)
we map to the equator of the sphere which corresponds to a harmonic map.
Hence, the ansatz \eqref{eq:definition-q-biharmonic} asks how far we need to "go away" from
the equator such that a harmonic map becomes biharmonic or conformal-biharmonic.

In particular, we want to study the following objectives:
\begin{enumerate}
    \item For which values of \(\alpha\) does the ansatz \eqref{eq:definition-q-biharmonic} lead to a biharmonic or a conformal-biharmonic map?
      \item Do we get any restriction on the harmonic map \(v\)?
    \item Are the biharmonic and conformal-biharmonic maps constructed via this approach stable?
\end{enumerate}

In the case of a closed domain and biharmonic maps, i.e. solutions of \eqref{eq:bitension}, the above questions can all be answered by the following

\begin{Theorem}
\label{thm:biharmonic-main}
Let \((M,g)\) be a closed Riemannian manifold and \(v\colon M\to\s^{n-1}\subset\R^{n}\) be a non-trivial harmonic map. 
Then, the map \(q\colon M\to\s^n\subset\R^{n+1}\) defined in \eqref{eq:definition-q-biharmonic}
is proper biharmonic if and only if 
\(\alpha=\frac{\pi}{4}\) and \(|\nabla v|^2=const\). 
In particular, the only proper biharmonic map obtained via the ansatz \eqref{eq:definition-q-biharmonic} is given by
\begin{align}
\label{eq:intro-thm-biharmonic-u}
q:=(\frac{1}{\sqrt{2}}\cdot v,  \frac{1}{\sqrt{2}}).  
\end{align}
The proper biharmonic maps constructed this way are unstable critical points of the bienergy \eqref{eq:bienergy}.
\end{Theorem}

As a second main result we will prove a similar theorem for conformal-biharmonic maps
between spheres.

\begin{Theorem}
\label{thm:c-biharmonic-main}
Let \((\s^m,g_{can})\) be the Euclidean sphere with the standard metric 
and consider a map of the form
\eqref{eq:definition-q-biharmonic} where \(v\colon \s^m\to\s^{n-1}\subset\R^{n}\) is a harmonic map.   Then, the map \(q\) is conformal-biharmonic if
\begin{align*}
|\nabla v|^2=\lambda,\qquad \sin^2\alpha
=\frac{1}{3}\frac{(m-1)(m-3)}{\lambda}+\frac{1}{2},    
\end{align*}
where \(\lambda\in\R\) is a constant.
The conformal-biharmonic maps constructed in this way are unstable critical points of the conformal bienergy \eqref{eq:bienergy-conformal}.
\end{Theorem}

In addition, motivated by the structure of biharmonic curves on spheres, we will perform the same analysis for maps
\(w\colon M\to\s^{n}\subset\R^{n+1}\)
of the form
\begin{align}
\label{eq:intro-ansatz-product}
w:=(\sin\beta\cdot v_1,\cos\beta\cdot v_2),\qquad \beta\in (0,\frac{\pi}{2}),  
\end{align}
where \(v_i\colon M\to\s^{n_i}\subset\R^{n_i+1},i=1,2\) are two harmonic maps and \(n_1+n_2=n-1\).
We will again investigate when maps of the form \eqref{eq:intro-ansatz-product}
can actually be biharmonic, i.e. a solution of \eqref{eq:bitension} or conformal-biharmonic, i.e. a solution of \eqref{eq-c-biharmonic-intro}.

In the biharmonic case we obtain the following result:

\begin{Theorem}
\label{thm:biharmonic-product-main}
Let \((M,g)\) be a closed Riemannian manifold and consider a map of the form \eqref{eq:intro-ansatz-product} where \(v_i\colon M\to\s^{n_i}\subset\R^{n_i+1},i=1,2\) are two harmonic maps and \(n_1+n_2=n-1\). Then, the map 
\(w\colon M\to\s^{n}\subset\R^{n+1}\) is proper biharmonic if and only if 
\begin{align*}
|\nabla v_1|^2-|\nabla v_2|^2=const\neq 0, 
 \qquad \beta=\frac{\pi}{4}  
\end{align*}
such that
\begin{align*}
w:=(\frac{1}{\sqrt{2}}\cdot v_1,  \frac{1}{\sqrt{2}}\cdot v_2).  
\end{align*}
The proper biharmonic maps constructed this way are unstable critical points of the bienergy \eqref{eq:bienergy}.
\end{Theorem}

Again, we will also provide a corresponding result for conformal-biharmonic maps
between Euclidean spheres.

\begin{Theorem}
\label{thm:c-biharmonic-product-main}
Let \((\s^m,g_{can})\) be the Euclidean sphere with the standard metric 
and consider a map of the form \eqref{eq:intro-ansatz-product} where \(v_i\colon\s^m\to\s^{n_i}\subset\R^{n_i+1},i=1,2\) are two harmonic maps
and \(n_1+n_2=n-1\).
Then, the map \(w\) is conformal-biharmonic if  
\begin{align}
|\nabla v_i|^2=\lambda_i=const,i=1,2,~~\lambda_1\neq\lambda_2,\qquad 
\cos 2\beta=\frac{2}{3}(m-1)(m-3)\frac{1}{\lambda_2-\lambda_1}. 
\end{align}   
\end{Theorem}

Let us make the following remarks on the main results obtained in this manuscript.

\begin{Bem}
\begin{enumerate}
    \item A result similar to Theorem \ref{thm:biharmonic-main} was obtained by 
    Loubeau and Oniciuc \cite[Theorem 1.1]{MR2327029}.
    More precisely, they considered a
     closed Riemannian manifold \((M,g)\) and a non-constant map \(\phi\colon M\to\s^n(\frac{r}{\sqrt{2}})\).
     They proved that the composite map \(u=\iota\circ\phi\colon M\to\s^{n+1}(r)\) is proper biharmonic if and only if \(\phi\) is harmonic and has constant energy density.
\item A weaker version of Theorem \ref{thm:biharmonic-main} was recently obtained by Ambrosie in \cite[Proposition 3.6]{MR4695264}
under the assumption that the harmonic map \(v\) has constant energy density.
    \item The result achieved in Theorem \ref{thm:biharmonic-main}  shows that in the case of a closed domain
     the only biharmonic map coming out of the ansatz 
     \eqref{eq:definition-q-biharmonic} 
     is the one discussed by Ou in \cite{MR2888980}. 
     More precisely, in \cite[Theorem 2.9]{MR2888980} Ou showed that the map defined in \eqref{eq:intro-thm-biharmonic-u} is proper biharmonic. Moreover, in 
     \cite[Remark 5]{MR2888980} he raised the question if there are other classes of biharmonic maps to spheres and the results provided by Theorem \ref{thm:biharmonic-main}  suggest that this is not the case.

    \item Recently, in \cite{MR4455953} Ou investigated the relation between 
bi-eigenmaps, buckling eigenmaps and biharmonic maps to spheres.
He assumes that there exists a biharmonic map \(\phi\colon M\to\s^n\) with constant energy density \(|d\phi|^2=c\) and, denoting by \(\iota\colon\s^n\to\R^{n+1}\) the canonical inclusion, that there also exists a function \(\rho\in C^\infty(M)\) satisfying \(\Delta^2(\iota\circ\phi)=-\rho\Delta(\iota\circ\phi)\).
Under these assumptions he showed (\cite[Theorem 3.3]{MR4455953}) that
the biharmonic map \(\phi\) is either harmonic or \(\rho=2c\).

    \item Biharmonic homogeneous polynomial maps between spheres have recently been constructed in \cite{MR4682222} and \cite{MR4599666}. The approach chosen by the authors has some similarity with the one used in this manuscript: They start with homogeneous polynomial maps, whose harmonicity is well-understood, and render them biharmonic in a second step.
     
     \item In the case that the domain \(M\) is non-compact or 
     has a boundary, the biharmonic map equation offers more flexibility
      as is demonstrated by the recent examples of proper biharmonic maps from the unit ball \cite{MR4830603}, \cite{MR4076824}.
     The reason for this difference is the following: In the case of a closed domain the maximum principle puts strong restrictions on the admissible values of \(\alpha\) in the ansätze \eqref{eq:definition-q-biharmonic}, \eqref{eq:intro-ansatz-product} and these restrictions are not present if the domain has a boundary,
     see the discussion at the end of Subsection \ref{subsec:w-biharmonic} for more details.

    \item There is another variant of biharmonic maps, the so-called \emph{extrinsic biharmonic maps}, which are favored in many analytic studies of the subject.
    As extrinsic biharmonic maps are obtained from a coercive energy functional 
    their existence can be guaranteed by variational techniques. On the other hand, it seems that the method of construction for biharmonic and conformal-biharmonic maps employed in this paper is not applicable to extrinsic biharmonic maps.
     For more details on extrinsic biharmonic maps we refer to the recent article \cite{MR4279391} and references therein.
       
    \item The statement of Theorem \ref{thm:biharmonic-main} on biharmonic maps
is more restrictive than the corresponding one on conformal-biharmonic maps, i.e. Theorem \ref{thm:c-biharmonic-main}. In the case of biharmonic maps all solutions constructed via the ansatz \eqref{eq:definition-q-biharmonic} need to have constant energy density while in the case of conformal-biharmonic maps there also might be solutions with non-constant energy density. This fact suggests that conformal-biharmonic maps are more flexible compared
to standard biharmonic maps.
\end{enumerate}
\end{Bem}

Throughout this article we will employ the following sign conventions: 
for the Riemannian curvature tensor field we use 
$$
R(X,Y)Z=[\nabla_X,\nabla_Y]Z-\nabla_{[X,Y]}Z,
$$ 
where \(X,Y,Z\) are vector fields.
Moreover, for the Ricci tensor field we employ
$$
\langle \Ric(X),Y\rangle=\Ric(X,Y)=\tr \left\{Z\to R(Z,X)Y\right\},
$$
and the scalar curvature is given by
$$
\Scal=\tr\Ric.
$$
For the rough Laplacian on the pull-back bundle $\phi^{\ast} TN$ we employ the analysts sign convention, i.e.
$$
\bar\Delta = \tr(\bar\nabla\bar\nabla-\bar\nabla_\nabla).
$$
In particular, this implies that the Laplace operator has a negative spectrum.

This article is organized as follows. In Section 2 we recall background material
on harmonic and biharmonic maps to spheres and derive the Euler-Lagrange equation
for conformal-biharmonic maps for a spherical target. In Section 3 we discuss
how the maps \eqref{eq:definition-q-biharmonic}, \eqref{eq:intro-ansatz-product}
can be turned into proper biharmonic maps and, finally, in Section 4 we address the same
question for conformal-biharmonic maps between spheres.

\medskip

\textbf{Acknowledgement}: The author would like to thank Cezar Oniciuc for many
inspiring comments on the manuscript and the unknown reviewers
for their many helpful suggestions.

\section{Preliminaries}
In this section we will recall the necessary background material that will be employed within this manuscript.

\subsection{Harmonic, biharmonic and conformal-biharmonic maps to spheres}
First, we recall the Euler-Lagrange equations for harmonic and biharmonic maps to spheres and, in addition, we will derive the Euler-Lagrange equation for conformal-biharmonic maps in the case of a spherical target.

Let \((M,g)\) be a Riemannian manifold and \(u\colon M\to\s^n\).
We consider the canonical embedding \(\iota\colon\s^n\hookrightarrow\R^{n+1}\) and will still write
\(u=(u_1,\ldots,u_{n+1})\) for \(\iota\circ u\).
In the following, we will employ the notation
\begin{align*}
\nabla u=(\nabla u_1,\ldots,\nabla u_{n+1}),
\qquad \Delta u=(\Delta u_1,\ldots,\Delta u_{n+1}).
\end{align*}
Here, \(\nabla\) refers to the gradient of \((M,g)\) and \(\Delta\)
is the Laplace-Beltrami operator. Note that each component of \(\nabla u\) is an \(m\)-dimensional vector, where \(m=\dim M\).

In this setup the energy of the map \(u\colon M\to\R^{n+1}\) acquires the form
\begin{align*}
E(u)=\frac{1}{2}\int_M|\nabla u|^2\dv,
\end{align*}
where we also impose the constraint \(|u|^2=1\), i.e. \(u\) takes values in \(\s^n\),
and its critical points are given by
\begin{align}
\label{eq:harmonic-sphere}
\Delta u+|\nabla u|^2u=0.    
\end{align}
Note that the second term, which represents the nonlinearity of the harmonic map equation, originates from orthogonally projecting the equation to the sphere and underlines the geometric nature of harmonic maps.

Recall that the bienergy functional for maps \(u\colon M\to\s^n\subset\R^{n+1}\) is given by
\begin{align}
\label{eq:bienergy-extrinsic}
E_2(u)=\frac{1}{2}\int_M(|\Delta u|^2-|\nabla u|^4)\dv,    
\end{align}
where we made use of the fact that
\begin{align*}
|\tau(u)|^2=\big|\Delta u+|\nabla u|^2u \big|^2
=|\Delta u|^2-|\nabla u|^4.
\end{align*}

The critical points of \eqref{eq:bienergy-extrinsic} are those which satisfy
\begin{align}
\label{eq:biharmonic-sphere}
\Delta^2u+2\operatorname{div}\big(|\nabla u|^2\nabla u\big)
-\big(\langle\Delta^2u,u\rangle-2|\nabla u|^4\big)u=0
\end{align}
and again the nonlinear terms proportional to \(u\) stem from projecting the equation onto the sphere. One could of course also derive \eqref{eq:biharmonic-sphere} by using the general biharmonic map equation \eqref{eq:bitension}
and inserting the geometric data of a spherical target.

Throughout this manuscript we will often employ the following statement.
\begin{Lem}
Let \(v\colon M\to\s^n\subset\R^{n+1}\) be a harmonic map.    
Then, the following formula holds
\begin{align}
\label{eq:v-harmonic-laplace}
0=&\Delta^2v+\Delta|\nabla v|^2v+2\nabla|\nabla v|^2\nabla v
-|\nabla v|^4v.
\end{align}
\end{Lem}

\begin{proof}
This follows from applying the Laplace operator to the equation for harmonic maps
to spheres \eqref{eq:harmonic-sphere}.
\end{proof}

As a final contribution of this section, we will derive the equation for conformal-biharmonic maps in the case of a spherical target in the extrinsic picture. 

In terms of the map \(u\colon M\to\s^n\subset\R^{n+1}\)
the conformal bienergy acquires the form
\begin{align}
\label{eq:c-energy-sphere}
E^c_2(u)=\frac{1}{2}\int_M\big(|\Delta u|^2-|\nabla u|^4+\frac{2}{3}\Scal^M|\nabla u|^2
-2\tr\langle \nabla u(\Ric^M(\cdot)),\nabla u(\cdot)\rangle\big)\dv
\end{align}
and we will now derive its critical points.

\begin{Prop}
The critical points of \eqref{eq:c-energy-sphere} are those which satisfy
\begin{multline}
\label{eq:euler-lagrange-cbi}
 \Delta^2u+2\operatorname{div}(|\nabla u|^2\nabla u)
+\frac{1}{3}\nabla u(\nabla\Scal^M)-\frac{2}{3}\Scal^M\Delta u
+2\tr\nabla^2u(\Ric^M(\cdot),\cdot) \\   
-\big(\langle\Delta^2u,u\rangle-2|\nabla u|^4+\frac{2}{3}\Scal^M|\nabla u|^2
-2\tr\langle\nabla u(\Ric^M(\cdot)),\nabla u(\cdot)\rangle\big)u=0.
\end{multline}
\end{Prop}

\begin{proof}
 We consider a variation of \(u\), that is a map 
 \(u_t\colon (-\epsilon,\epsilon)\times M\to\s^n\subset\R^{n+1},\epsilon>0\) 
 satisfying \(d u_t(\partial_t)\big|_{t=0}=\eta\).
Moreover, we introduce the Lagrange multiplyer \(\Lambda\), which encodes that the image of \(u\) is constrained to the unit sphere, and consider
\begin{align*}
E^c_2(u)=\frac{1}{2}\int_M\big(|\Delta u|^2-|\nabla u|^4+\frac{2}{3}\Scal^M|\nabla u|^2
-2\tr\langle \nabla u(\Ric^M(\cdot)),\nabla u(\cdot)\rangle+\Lambda(|u|^2-1)\big)\dv.
\end{align*}

Now, a direct calculation shows that
\begin{align*}
\frac{d}{dt}\big|_{t=0}E^c_2(u_t)=
\int_M\big(&\langle\eta,\Delta^2u+2\operatorname{div}(|\nabla u|^2\nabla u)
+\frac{2}{3}\Scal^M\langle\nabla\eta,\nabla u\rangle \\
&-2\tr\langle \nabla u(\Ric^M(\cdot)),\nabla \eta(\cdot)\rangle
+\Lambda u\big)\dv.
\end{align*}
The variation with respect to the Lagrange multiplyer \(\Lambda\) yields the constraint \(|u|^2=1\).

Now, integration by parts gives
\begin{align*}
\int_M\Scal^M\langle\nabla\eta,\nabla u\rangle\dv&=-\int_M\big(\Scal^M\langle\eta,\Delta u\rangle
+\langle\eta,\nabla u(\nabla\Scal^M)\rangle\big)\dv,\\
\int_M\tr\langle \nabla u(\Ric^M(\cdot)),\nabla \eta(\cdot)\rangle\dv
&=-\int_M\big(\langle\eta,\tr\nabla^2u(\Ric^M(\cdot),\cdot)+\frac{1}{2}\nabla u(\nabla\Scal^M)\rangle\big)\dv,
\end{align*}
where we used the identity \(\operatorname{div}\Ric=\frac{1}{2}\nabla\Scal g\) which holds on every Riemannian manifold.
Combining these equations we find that the first variation leads to 
\begin{align*}
\Delta^2u+2\operatorname{div}(|\nabla u|^2\nabla u)
+\frac{1}{3}\nabla u(\nabla\Scal^M)-\frac{2}{3}\Scal^M\Delta u
+2\tr\nabla^2u(\Ric^M(\cdot),\cdot)+\Lambda u=0.
\end{align*}
In order to eliminate the Lagrange multiplyer \(\Lambda\) we test the above equation with \(u\)
and find
\begin{align*}
-\Lambda=\langle\Delta^2u,u\rangle-2|\nabla u|^4+\frac{2}{3}\Scal^M|\nabla u|^2
-2\tr\langle\nabla u(\Ric^M(\cdot)),\nabla u(\cdot)\rangle,
\end{align*}
reinserting the expression for \(\Lambda\) then completes the proof.
\end{proof}

Note that we could also derive the equation for conformal-biharmonic maps to spheres
by starting from the general equation for conformal-biharmonic maps \eqref{eq-c-biharmonic-intro} and using the geometric data of a spherical target.

\subsection{Eigenmaps to Euclidean spheres}
In this subsection we briefly recall a number of general results on eigenmaps, for more details 
on this subject we refer to the books of Eells \(\&\) Ratto \cite[Part 3]{MR1242555}
and Baird \(\&\) Wood \cite[p. 79]{MR2044031}.

\begin{Dfn}
A smooth map \(\phi\colon M\to\s^n\) is called eigenmap if the components
of \(u:=\iota\circ\phi\colon M\to\s^n\subset\R^{n+1}\) are all eigenfunctions of the Laplace-Beltrami operator on 
\(M\) with the same eigenvalue.
\end{Dfn}

An immediate consequence of the above definition is the following fact
\begin{Cor}
A smooth map \(\phi\colon M\to\s^n\) is an eigenmap if and only if it is harmonic
with constant energy density \(e(\phi):=\frac{1}{2}|d\phi|^2\).
\end{Cor}

\begin{Prop}
\label{prop:homo-pol}
Let \(F\colon\R^{m+1}\to\R^{n+1}\) be a harmonic map whose components are homogeneous
polynomials of the same degree \(k\in\{1,2,\ldots\}\).
Moreover, suppose that \(F\) restricts to a map \(\phi\colon\s^m\to\s^n\).
Then \(\phi\) is harmonic with constant energy density
\begin{align}
e(\phi)=\frac{k}{2}(k+m-1).
\end{align}
\end{Prop}

Fore more details on polynomial maps between spheres and how they can be 
used to construct biharmonic maps we refer to \cite{MR4599666}.

\section{From harmonic to biharmonic maps}
In this section we will discuss an algorithm that produces a biharmonic map 
from either one given harmonic map or two given harmonic maps.

\subsection{Rendering a harmonic map biharmonic}
\label{subsec:w-biharmonic}
First, we start with the case that we have a given harmonic map and aim at rendering it proper biharmonic.
More precisely, in order to construct a solution to \eqref{eq:biharmonic-sphere}
we consider the map \(q\colon M\to\s^{n}\subset\R^{n+1}\) which was defined in 
\eqref{eq:definition-q-biharmonic} as follows
\begin{align*}
q:=\big(\sin\alpha\cdot v,\cos\alpha\big).   
\end{align*}
Here, we assume that \(v\colon M\to\s^{n-1}\subset\R^{n}\) is a given harmonic map, i.e. a solution of \eqref{eq:harmonic-sphere},
and \(\alpha\in (0,\frac{\pi}{2})\) a real parameter.

The next Proposition gives a characterization when maps of the form \eqref{eq:definition-q-biharmonic} are actually biharmonic.

\begin{Prop}
The map \(q\colon M\to\s^n\subset\R^{n+1}\) defined in \eqref{eq:definition-q-biharmonic} is biharmonic if and only if the following equation holds
\begin{align}
\label{eq:biharmonic-w}
(\Delta|\nabla v|^2)v=-2\nabla(|\nabla v|^2)\nabla v
+2\sin^2\alpha\operatorname{div}(|\nabla v|^2\nabla v)+|\nabla v|^4v.
\end{align}
\end{Prop}

\begin{proof}
Inserting \eqref{eq:definition-q-biharmonic} into the equation for biharmonic maps to spheres \eqref{eq:biharmonic-sphere} we get the following constraint from the last component
\begin{align*}
\big(\langle\Delta^2v,v\rangle-2\sin^2\alpha|\nabla v|^4\big)\sin^2\alpha\cos\alpha=0.    
\end{align*}
Suppose that this constraint holds, the first \(n\)-components yield
\begin{align*}
\sin\alpha\big(\Delta^2v+2\sin^2\alpha\nabla(|\nabla v|^2\nabla v)\big)=0.
\end{align*}
Using the identity 
\begin{align*}
\langle \operatorname{div}(|\nabla v|^2\nabla v),v\rangle    
=-|\nabla v|^4,
\end{align*}
which holds for every map \(v\colon M\to\s^{n-1}\subset\R^{n}\),
it can easily be seen that the equation for the \((n+1)\)-th component of \(q\)
follows from the equation for the first \(n\)-components and thus holds automatically if \eqref{eq:biharmonic-w} is satisfied.

By assumption \(v\colon M\to\s^{n-1}\subset\R^{n}\) is a harmonic map, i.e. a solution of \eqref{eq:harmonic-sphere}, such that we can apply \eqref{eq:v-harmonic-laplace}.
Combining this identity with the previous equation now completes the proof.
\end{proof}

Note that equation \eqref{eq:biharmonic-w} is only of third order whereas the biharmonic map equation \eqref{eq:biharmonic-sphere} is of fourth order. Hence, by starting with a harmonic map with the aim of deforming it to a biharmonic one we already achieved to reduce the order of derivatives by one.

Later on, we will see that we can further reduce the number of derivatives in \eqref{eq:biharmonic-w} but it turns out that \eqref{eq:biharmonic-w} is particularly
suited in order to detect the restrictions that arise when deforming a harmonic map to a biharmonic one via the ansatz \eqref{eq:biharmonic-w}.

\begin{Prop}
Let \((M,g)\) be a Riemannian manifold
and \(v\colon M\to\s^{n-1}\subset \R^{n}\) be a solution of \eqref{eq:biharmonic-w}. Then, the following constraint must be satisfied
\begin{align}
\label{eq:constraint-w}
\Delta|\nabla v|^2=(1-2\sin^2\alpha)|\nabla v|^4.
\end{align}
\end{Prop}

\begin{proof}
The claim follows by testing \eqref{eq:biharmonic-w} with \(v\).
\end{proof}

A direct calculation now shows that, in the case of a closed domain, the only way to produce biharmonic maps via the ansatz \eqref{eq:definition-q-biharmonic} is covered by the following:

\begin{Theorem}
\label{thm:biharmonic-w}
Let \((M,g)\) be a closed Riemannian manifold and \(v\colon M\to\s^{n-1}\subset\R^{n}\) be a non-trivial harmonic map. 
Then, the map \(q\colon M\to\s^n\subset\R^{n+1}\) defined in \eqref{eq:definition-q-biharmonic}
is proper biharmonic if and only if 
\(\alpha=\frac{\pi}{4}\) and \(|\nabla v|^2=const\). 
In particular, the only proper biharmonic map obtained via the ansatz \eqref{eq:definition-q-biharmonic} is given by
\begin{align}
\label{eq:w-biharmonic-sqrt2}
q:=\big(\frac{1}{\sqrt{2}} \cdot v,\frac{1}{\sqrt{2}}\big).
\end{align}
\end{Theorem}

\begin{proof}
Since \(M\) is closed by assumption we can make use of the maximum principle, and we will now apply it to 
\eqref{eq:constraint-w} as follows:
If \(1-2\sin^2\alpha>0\) the maximum principle implies that \(|\nabla v|^2=const\)
and reinserting in \eqref{eq:constraint-w} then yields that \(|\nabla v|^2=0\) meaning that \(v\) is trivial which we have excluded in the assumptions.
If \(1-2\sin^2\alpha<0\) we can apply the maximum principle to \(-|\nabla v|^2\)
and will arrive at the same conclusion. In the remaining case that \(\alpha=\frac{\pi}{4}\)
the right hand side of \eqref{eq:constraint-w} vanishes and \(|\nabla v|^2\) has to be a harmonic function,
and as \(M\) is closed, we get that \(|\nabla v|^2=const\).
Finally, it can be checked directly that \(\alpha=\frac{\pi}{4}\) and \(|\nabla v|^2=const\)
provides a solution of \eqref{eq:biharmonic-w} completing the proof.
\end{proof}

The biharmonic maps constructed in \eqref{eq:w-biharmonic-sqrt2} have bienergy
\begin{align*}
E_2(w)=\frac{|\nabla v|^4}{8}\vol(M,g)=const.
\end{align*}
This strongly suggests that these maps are all unstable critical points of the bienergy \eqref{eq:bienergy}.

Indeed, we can apply a general instability result for biharmonic maps
to show that the biharmonic maps provided by Theorem \ref{thm:biharmonic-w} are actually
unstable critical points of the bienergy \eqref{eq:bienergy}. 

\begin{Theorem}
The proper biharmonic maps obtained in Theorem \ref{thm:biharmonic-w}
are unstable critical points of the bienergy \eqref{eq:bienergy}.
\end{Theorem}

\begin{proof}
We apply the following well-established result for biharmonic maps:
In \cite[Theorem 2.5]{MR4110268}, see also \cite[Theorem 16]{MR886529}, it was
shown that a biharmonic map from a closed Riemannian manifold into a space form of positive sectional curvature with \(|\tau(\phi)|^2=const\) must be unstable.
It can easily be checked that the norm of the tension field of the proper biharmonic maps provided by Theorem \ref{thm:biharmonic-w} is constant such that the claim follows from the aforementioned general result.
\end{proof}

This completes the proof of Theorem \ref{thm:biharmonic-main}.

\begin{Bsp}
In the case that \(\phi\colon\s^n(r)\to\s^{n+1}\) is an isometric immersion,
we have \(|d\phi|^2=n\) and Theorem \ref{thm:biharmonic-w} provides the small proper biharmonic hypersphere \(\phi\colon\s^n(\frac{1}{\sqrt{2}})\to\s^{n+1}\).
\end{Bsp}

The next Proposition shows that we can further reduce \eqref{eq:biharmonic-w}
to a second order equation in \(v\) which, unfortunately, is not very useful for explicit calculations.

\begin{Prop}
The map \(q\colon M\to\s^n\subset\R^{n+1}\) defined in \eqref{eq:definition-q-biharmonic} is biharmonic if and only if
\begin{align}
\label{eq:cor-w-second}
|\nabla^2v|^2v=-\Ric^M(\nabla v,\nabla v)v-\cos^2\alpha \nabla(|\nabla v|^2)\nabla v
+\sin^2\alpha |\nabla v|^4v+\frac{1}{2}|\nabla v|^4v.
\end{align}
Here, \(\Ric^M\) is interpreted as a \((0,2)\)-tensor.
\end{Prop}
\begin{proof}
A direct calculation with the help of the Bochner formula shows that
\begin{align*}
\Delta|\nabla v|^2=&2|\nabla^2v|^2+2\Ric^M(\nabla v,\nabla v)+2\langle\nabla v,\nabla\Delta v\rangle \\
=&2|\nabla^2v|^2+2\Ric^M(\nabla v,\nabla v)-2|\nabla v|^4,
\end{align*}
where we used that \(v\) is a harmonic map, that is a solution of \eqref{eq:harmonic-sphere}.

In addition, a direct calculation shows that
\begin{align*}
\operatorname{div}(|\nabla v|^2\nabla v)=\nabla(|\nabla v|^2)\nabla v+|\nabla v|^2\Delta v=\nabla(|\nabla v|^2)\nabla v-|\nabla v|^4v.
\end{align*}
The claim now follows from combining all the above identities.
\end{proof}

\begin{Prop}
Let \((M,g)\) be a closed Riemannian manifold. Then \eqref{eq:cor-w-second} implies
\begin{align*}
\Ric^M(\nabla v,\nabla v)\leq \frac{m-1}{m}|\nabla v|^4.    
\end{align*}
In particular, if \(M=\s^m\) with the round metric then \(m\leq|\nabla v|^2\).
\end{Prop}

\begin{proof}
Testing \eqref{eq:cor-w-second} with \(v\) and choosing \(\alpha=\frac{\pi}{4}\),    
which needs to hold in the case of a closed domain
due to Theorem \ref{thm:biharmonic-w}, we find
\begin{align*}
|\nabla^2v|^2v=-\Ric^M(\nabla v,\nabla v)v+|\nabla v|^4v.    
\end{align*}
Now, from the Cauchy-Schwarz inequality and the fact that \(v\) is a harmonic map
we get
\begin{align*}
|\nabla v|^4=|\Delta v|^2\leq m|\nabla^2v|^2    
\end{align*}
and combining both equations completes the proof.
\end{proof}

As a last important observation in this section we want to point out that
for \(M\) being non-compact or compact with boundary the maximum principle
no longer imposes the constraints that led to the statement of Theorem 
\ref{thm:biharmonic-w} such that the ansatz \eqref{eq:definition-q-biharmonic} offers more flexibility in this case.

Note that by choosing \(\alpha=\frac{\pi}{4}\) the constraint 
\eqref{eq:constraint-w} yields that \(|\nabla v|^2\) needs to be a harmonic function
and for \(M\) being compact this forces \(|\nabla v|^2\) to be constant.
However, in the non-compact case harmonic functions are not constraint by the maximum principle such that we can directly construct additional biharmonic maps out of the ansatz \eqref{eq:definition-q-biharmonic}.

To this end, we consider the following three maps given by
\begin{align}
\label{eq:three-harmonic-maps}
\pi&=\big(\pi_i:=\frac{x_i}{r}\big)_{1\leq i\leq m}, \\
\nonumber \mu&=\bigg(\mu_{ij}:=\frac{1}{\sqrt{m(m-1)}}\big(-\delta_{ij}+m\frac{x_ix_j}{r^2}\big)\bigg)_{1\leq i,j\leq m},\\
\nonumber \nu&=\bigg(\nu_{ijk}:=\frac{1}{\sqrt{(m-1)(m+2)}}
\big(\delta_{ij}\frac{x_k}{r}+\delta_{jk}\frac{x_i}{r}
+\delta_{ik}\frac{x_j}{r}-(m+2)\frac{x_ix_jx_k}{r^3}\big)\bigg)_{1\leq i,j,k\leq m},
\end{align}
where 
\begin{align*}
\pi\colon\R^m\setminus\{0\}\to\s^{m-1},\qquad 
\mu\colon\R^m\setminus\{0\}\to\s^{m^2-1},\qquad
\nu\colon\R^m\setminus\{0\}\to\s^{m^3-1}.
\end{align*}
Here, \(r:=\sqrt{x_1^2+\ldots +x_m^2}\) represents the Euclidean distance.
These comprise harmonic maps to the sphere, see \cite{MR4371934} for the details, and can be rendered biharmonic via the ansatz \eqref{eq:definition-q-biharmonic}, 
see \cite{{MR4830603}}, \cite{MR4076824} for the explicit calculations,  
which we summarize within the following list:

\begin{enumerate}
\item The map \(\tilde \pi\colon\R^m\setminus\{0\}\to\s^{m}\subset\R^{m+1}\) defined by
\begin{align*}
\tilde \pi:=(\sin\alpha\cdot\pi,\cos\alpha),\qquad \alpha\in (0,\frac{\pi}{2})
\end{align*}
  is proper biharmonic if and only if
\begin{align*}
\sin^2\alpha=\frac{3}{2}\frac{m-3}{m-1},\qquad 4\leq m\leq 6.   
\end{align*}
It can also be considered as a weak biharmonic map 
from \(B^m\to\s^m,m=5,6\), see \cite[Theorem 1.1]{MR4076824}
where \(B^m\) represents the unit ball in \(m\)-dimensions.
\item The map \(\tilde \mu\colon\R^m\setminus\{0\}\to\s^{m^2}\subset\R^{m^2+1}\) defined by
\begin{align*}
\tilde \mu:=(\sin\beta\cdot \mu,\cos\beta)    ,\qquad \beta\in (0,\frac{\pi}{2})
\end{align*}
is proper biharmonic if and only if 
\begin{align*}
\sin^2\beta=\frac{m-2}{m},\qquad m\geq 3.
\end{align*}
It can also be thought of as a weak biharmonic map 
from \(B^m,m\geq 5\), see \cite[Theorem 1.2]{MR4830603}.
\item The map \(\tilde \nu\colon\R^m\setminus\{0\}\to\s^{m^3}\subset\R^{m^3+1}\) defined by
\begin{align*}
\tilde \nu:=(\sin\gamma\cdot \nu,\cos\gamma),\qquad \gamma\in (0,\frac{\pi}{2}) \end{align*}
is proper biharmonic if and only if 
\begin{align*}
\sin^2\gamma=\frac{5}{6}\frac{m-1}{m+1},\qquad m\geq 2.
\end{align*}
It can also be interpreted as a weak biharmonic map 
from \(B^m,m\geq 5\), see \cite[Theorem 1.4]{MR4830603}.
\end{enumerate}

It should be noted that in four dimensions all of the above proper biharmonic maps
have a particular structure which is very similar to the case of biharmonic maps from a closed domain, in particular the last component always equals \(\frac{1}{\sqrt{2}}\).
More precisely, for \(m=4\) we find that the proper biharmonic maps
from the above list are of the form
\begin{align*}
\tilde\pi&:=\frac{1}{\sqrt{2}}\bigg(\big(\frac{x_i}{r}\big)_{1\leq i\leq m},1\bigg),\\
\tilde \mu&:=\frac{1}{\sqrt{2}}\bigg(\frac{1}{\sqrt{12}}\big(-\delta_{ij}+4\frac{x_ix_j}{r^2}\big)_{1\leq i,j\leq m},1\bigg),\\
\tilde \nu&:=\frac{1}{\sqrt{2}}\bigg(
\frac{1}{\sqrt{18}}
\big(\delta_{ij}\frac{x_k}{r}+\delta_{jk}\frac{x_i}{r}
+\delta_{ik}\frac{x_j}{r}-6\frac{x_ix_jx_k}{r^3}\big)_{1\leq i,j,k\leq m},1\bigg).
\end{align*}
Here, we consider all these maps as maps from \(\R^4\setminus\{0\}\)
and it is evident that their structure is similar to the compact case.
Indeed, these maps arise if we choose \(\alpha=\frac{\pi}{4}\) in the ansatz
\eqref{eq:definition-q-biharmonic} but now \(|\nabla\tilde p|^2, |\nabla\tilde\mu|^2, |\nabla\tilde\nu|^2\) are non-constant harmonic functions on the domain. Although the biharmonic map equation seems to offer more flexibility in the case of a non-compact domain it is quite remarkable that 
in four dimensions, which is the critical dimension, the structure of the explicit solutions in the compact and the non-compact case has some common features, namely
the appearance of the overall factor of \(\frac{1}{\sqrt{2}}\).

\subsection{Manufacturing a biharmonic out of two harmonic maps}
As a next step we now make an ansatz motivated from the second 
class of proper biharmonic curves presented in the introduction.

More precisely, we investigate when maps of the form
\(w\colon M\to\s^n\subset\R^{n+1}\) (see \eqref{eq:intro-ansatz-product} in the introduction) defined by 
\begin{align*}
w:=(\sin\beta\cdot v_1,\cos\beta\cdot v_2),\qquad \beta\in(0,\frac{\pi}{2}),
\end{align*}
where \(v_i\colon M\to\s^{n_i}\subset\R^{n_i+1},i=1,2,n_1+n_2=n-1\)
are two given harmonic maps,
are proper biharmonic.
In particular, we want to determine conditions on 
\(\beta\) and \(v_i,i=1,2\) that render the map \eqref{eq:intro-ansatz-product}
proper biharmonic.

\begin{Prop}
The map \(\eqref{eq:intro-ansatz-product}\) is biharmonic if and only if
\begin{align}
\label{eq:equation-biharmonic-product}
\big(\Delta|\nabla v_i|^2\big)v_i=&-2(\nabla|\nabla v_i|^2)\nabla v_i    
+|\nabla v_i|^4v_i
+2\nabla(\sin^2\beta|\nabla v_1|^2+\cos^2\beta|\nabla v_2|^2)\nabla v_i \\
\nonumber &-2(\sin^2\beta|\nabla v_1|^2+\cos^2\beta|\nabla v_2|^2)|\nabla v_i|^2v_i \\
\nonumber &+\big(\sin^2\beta(\Delta|\nabla v_1|^2-|\nabla v_1|^4)
+\cos^2\beta(\Delta|\nabla v_2|^2-|\nabla v_2|^4) \\
\nonumber &+2\sin^4\beta|\nabla v_1|^4+2\cos^4\beta|\nabla v_2|^4
+4\sin^2\beta\cos^2\beta|\nabla v_1|^2|\nabla v_2|^2\big)v_i,
\end{align}
where \(i=1,2\).
\end{Prop}

\begin{proof}
By definition \eqref{eq:intro-ansatz-product} we immediately get
\begin{align*}
|\nabla w|^2=&\sin^2\beta|\nabla v_1|^2+\cos^2\beta|\nabla v_2|^2,\\
\langle\Delta^2w,w\rangle=&\sin^2\beta\langle\Delta^2v_1,v_1\rangle
+\cos^2\beta\langle\Delta^2v_2,v_2\rangle.
\end{align*}
The result now follows by a direct calculation using \eqref{eq:v-harmonic-laplace}
which can be applied as \(v_i,i=1,2\) are both harmonic by assumption.
\end{proof}

The next Proposition shows that the system
\eqref{eq:equation-biharmonic-product}
leads to the following restrictions on \(\beta\)
and the harmonic maps \(v_i\colon M\to\s^{n_i}\subset\R^{n_i+1},i=1,2\).

\begin{Prop}
\label{prop:product-biharmonic-constraint}
Let \((M,g)\) be a Riemannian manifold and  
consider a map \(w\colon M\to\s^n\subset\R^{n+1}\) as defined in 
\eqref{eq:intro-ansatz-product}. Then, in order to obtain a biharmonic map,
the following constraint needs to be satisfied
\begin{align}
\label{eq:constraints-biharmonic-product}
\Delta\big(|\nabla v_1|^2-|\nabla v_2|^2\big)
&=\cos 2\beta\big(|\nabla v_1|^2-|\nabla v_2|^2\big)^2.
\end{align}
\end{Prop}

\begin{proof}
Testing \eqref{eq:equation-biharmonic-product} with \(v_i,i=1,2\) we get the two equations
\begin{align*}
\cos^2\beta\big(\Delta|\nabla v_1|^2-\Delta|\nabla v_2|^2\big)=&
(1-3\sin^2\beta+2\sin^4\beta)|\nabla v_1|^4
+\cos^2\beta(2\cos^2\beta-1)|\nabla v_2|^4 \\
\nonumber&+2\cos^2\beta(2\sin^2\beta-1)|\nabla v_1|^2|\nabla v_2|^2,\\
\nonumber\sin^2\beta\big(\Delta|\nabla v_2|^2-\Delta|\nabla v_1|^2\big)=&
(1-3\cos^2\beta+2\cos^4\beta)|\nabla v_2|^4
+\sin^2\beta(2\sin^2\beta-1)|\nabla v_1|^4 \\
\nonumber&+2\sin^2\beta(2\cos^2\beta-1)|\nabla v_1|^2|\nabla v_2|^2.
\end{align*}
Now, using the trigonometric identities
\begin{align*}
1-3\sin^2\beta+2\sin^4\beta=&\cos^2\beta\cos 2\beta, \\
1-3\cos^2\beta+2\cos^4\beta=&-\sin^2\beta\cos 2\beta
\end{align*}
we obtain the result.
\end{proof}

In the case that the domain is a closed Riemannian manifold we obtain the following 

\begin{Theorem}
\label{thm:biharmonic-product}
Let \((M,g)\) be a closed Riemannian manifold and consider a map of the form \eqref{eq:intro-ansatz-product} where \(v_i\colon M\to\s^{n_i}\subset\R^{n_i+1},i=1,2\) are two harmonic maps and 
\(n_1+n_2=n-1\). Then, the map 
\(w\colon M\to\s^{n}\subset\R^{n+1}\) is proper biharmonic if and only if 
\begin{align*}
|\nabla v_1|^2-|\nabla v_2|^2=const\neq 0, 
 \qquad \beta=\frac{\pi}{4}  
\end{align*}
such that
\begin{align}
\label{eq:biharmonic-product-solution}
w:=(\frac{1}{\sqrt{2}}\cdot v_1,  \frac{1}{\sqrt{2}}\cdot v_2).  
\end{align}
\end{Theorem}

\begin{proof}
Integrating \eqref{eq:constraints-biharmonic-product} over \(M\) and using the 
divergence theorem we find
\begin{align*}
0=\cos 2\beta\int_M\big(|\nabla v_1|^2-|\nabla v_2|^2\big)^2\dv.    
\end{align*}
If \(|\nabla v_1|^2=|\nabla v_2|^2\) the above constraint surely holds true.
On the other hand, if \(\beta=\frac{\pi}{4}\), then the above equation is satisfied as well and applying the maximum principle to \eqref{eq:constraints-biharmonic-product} then yields that \(|\nabla v_1|^2-|\nabla v_2|^2=const\).

We will show that the first case does not occur. If it would hold, then we have \(|\nabla w|^2=|\nabla v_1|^2=|\nabla v_2|^2\). Now, a direct calculation shows that
\begin{align*}
\Delta w=(\sin\beta\cdot\Delta v_1,\cos\beta\cdot\Delta v_2)
=-(\sin\beta\cdot|\nabla v_1|^2v_1,\cos\beta\cdot
|\nabla v_2|^2v_2)=-|\nabla w|^2w    
\end{align*}
which means that \(w\) is actually harmonic such that we do not obtain a proper biharmonic map.

Regarding the second case, inserting \(\beta=\frac{\pi}{4}\) into  
\eqref{eq:equation-biharmonic-product}, we find
\begin{align*}
\frac{1}{2}\Delta\big(|\nabla v_1|^2-|\nabla v_2|^2\big)v_i
=-\nabla\big(|\nabla v_1|^2-|\nabla v_2|^2\big)\nabla v_i,\qquad i=1,2
\end{align*}
and it is clear that \(|\nabla v_1|^2-|\nabla v_2|^2=const\)  
provides a solution to this equation.
\end{proof}

\begin{Bem}
Note that we do not require that both \(|\nabla v_1|^2,|\nabla v_2|^2\) are constant
in the previous Theorem, it is only their difference which needs to be constant.
\end{Bem}

Finally, we note that the norm of the tension field of the
proper biharmonic map \eqref{eq:biharmonic-product-solution} is equal to a constant
and thus these biharmonic maps are again unstable critical points of the bienergy \eqref{eq:bienergy}.

This completes the proof of Theorem \ref{thm:biharmonic-product-main}.
 
Again, it is advantageous to consider a domain that is non-compact when it comes to constructing proper biharmonic maps
via the ansatz \eqref{eq:intro-ansatz-product} in order to circumvent 
restrictions due to the maximum principle. 
More precisely, in the case of a closed
domain all solutions constructed by the ansatz \eqref{eq:intro-ansatz-product} are covered by Theorem \ref{thm:biharmonic-product} but there are additional solutions in the case of a non-compact domain as is evidenced by the following Theorem.

\begin{Theorem}
Consider the three harmonic maps  \(\pi,\mu,\nu\) detailed in \eqref{eq:three-harmonic-maps}. Then, the map
\begin{enumerate}
    \item \(w_1\colon\R^m\setminus\{0\}\to\s^{m+m^2-2}\) defined by
    \begin{align*}
    w_1:=(\sin\alpha\cdot\pi,\cos\alpha\cdot\mu),\qquad \alpha\in(0,\frac{\pi}{2})
\end{align*}
is proper biharmonic if
\begin{align*}
 \cos 2\alpha=2\frac{m-4}{m+1}   
\end{align*}
which implies that \(3<m<9\).
\medskip
    \item \(w_2\colon\R^m\setminus\{0\}\to\s^{m+m^3-2}\) defined by
    \begin{align*}
    w_2:=(\sin\beta\cdot\pi,\cos\beta\cdot\nu),\qquad \beta\in(0,\frac{\pi}{2})
\end{align*}
is proper biharmonic if
\begin{align*}
 \cos 2\beta=\frac{m-4}{m+2}
\end{align*}
which can be solved for any \(m>1\).
\medskip
\item \(w_3\colon\R^m\setminus\{0\}\to\s^{m^2+m^3-2}\) defined by
    \begin{align*}
    w_3:=(\sin\gamma\cdot \mu,\cos\gamma\cdot \nu),\qquad \gamma\in(0,\frac{\pi}{2})
\end{align*}
is proper biharmonic if
\begin{align*}
 \cos 2\gamma=2\frac{m-4}{m+3}
\end{align*}
which implies that \(2<m<11\).
\end{enumerate}
\end{Theorem}

\begin{proof}
First of all, we recall that
\begin{align*}
|\nabla\pi|^2=\frac{m-1}{r^2},\qquad |\nabla\mu|^2=\frac{2m}{r^2},\qquad   
|\nabla\nu|^2=3\frac{m+1}{r^2}. 
\end{align*}
Moreover, we have the identity
\begin{align*}
  \Delta^{\R^m}\frac{1}{r^2}=2\frac{4-m}{r^4},
\end{align*}
where \(\Delta^{\R^m}\) represents the Laplace operator on \(\R^m\).

Inserting these identities into the constraint \eqref{eq:constraints-biharmonic-product} then yields the equations for \(\alpha,\beta,\gamma\).
Now, it remains to show that the maps \(w_j,j=1,2,3\) are actual solutions of the biharmonic map equation \eqref{eq:equation-biharmonic-product}.
However, it is straightforward to see that
\begin{align*}
 \langle\nabla|\nabla v_j|^2,\nabla v_k\rangle=0,\qquad j,k=1,2,3,   
\end{align*}
where we denote \(v_1=\pi,v_2=\mu,v_3=\nu\). Hence for the maps under consideration the equation for biharmonic maps \eqref{eq:equation-biharmonic-product} is satisfied whenever the constraint \eqref{eq:constraints-biharmonic-product} holds.
\end{proof}

\begin{Bem}
In the case of maps from \(\R^4\setminus\{0\}\) the three proper 
biharmonic maps obtained in the previous Theorem acquire the form
\begin{align*}
w_1&=\frac{1}{\sqrt{2}}\bigg(\big(\frac{x_i}{r}\big)_{1\leq i\leq m},\frac{1}{\sqrt{12}}\big(-\delta_{ij}+4\frac{x_ix_j}{r^2}\big)_{1\leq i,j\leq m}\bigg),\\
w_2&=\frac{1}{\sqrt{2}}\bigg(\big(\frac{x_i}{r}\big)_{1\leq i\leq m},\frac{1}{\sqrt{18}}
\big(\delta_{ij}\frac{x_k}{r}+\delta_{jk}\frac{x_i}{r}
+\delta_{ik}\frac{x_j}{r}-6\frac{x_ix_jx_k}{r^3}\big)_{1\leq i,j,k\leq m}\bigg),\\    
w_3&=\frac{1}{\sqrt{2}}\bigg(\frac{1}{\sqrt{12}}\big(-\delta_{ij}+4\frac{x_ix_j}{r^2}\big)_{1\leq i,j\leq m},\frac{1}{\sqrt{18}}
\big(\delta_{ij}\frac{x_k}{r}+\delta_{jk}\frac{x_i}{r}
+\delta_{ik}\frac{x_j}{r}-6\frac{x_ix_jx_k}{r^3}\big)_{1\leq i,j,k\leq m} \bigg).
\end{align*}
Hence, as already noticed at the end of Subsection \ref{subsec:w-biharmonic},
the well-known overall factor of \(\frac{1}{\sqrt{2}}\) again appears in the critical dimension.
\end{Bem}

\section{Constructing conformal-biharmonic between spheres}
In this section we will carry out a similar analysis as in the previous section, i.e.
we will investigate if we can deform a harmonic, or a set of two harmonic maps
to a solution of the equation for conformal-biharmonic maps to spheres \eqref{eq:euler-lagrange-cbi}.

We start by considering the map \(q\colon M\to\s^n\subset\R^{n+1}\), see \eqref{eq:definition-q-biharmonic} in the introduction, given by
\begin{align*}
q:=\big(\sin\alpha\cdot v,\cos\alpha\big).
\end{align*}
Here, we again assume that \(v\colon M\to\s^{n-1}\subset\R^{n}\) is a harmonic map, i.e. a solution of \eqref{eq:harmonic-sphere},
and \(\alpha\in (0,\frac{\pi}{2})\) a real parameter.

\begin{Prop}
The map \(q\colon M\to\s^n\subset\R^{n+1}\) defined in \eqref{eq:definition-q-biharmonic}
is conformally biharmonic if the following equation holds
\begin{align}
\label{eq:w-cbi}
(\Delta|\nabla v|^2)v=&-2\cos^2\alpha(\nabla|\nabla v|^2)\nabla v
+(1-2\sin^2\alpha)|\nabla v|^4v
+\frac{1}{3}\nabla v(\nabla\Scal^M) \\
\nonumber&+\frac{2}{3}\Scal^M|\nabla v|^2v
+2\tr\nabla^2v(\Ric^M(\cdot),(\cdot)).
\end{align}
\end{Prop}

\begin{proof}
This follows from \eqref{eq:euler-lagrange-cbi} and the fact that \(v\)
is a harmonic map.
\end{proof}

It is obvious that we do not get much information out of \eqref{eq:w-cbi} as long
as we do not have an explicit expression for the Ricci curvature on \(M\).
On the other hand \eqref{eq:w-cbi} allows us to obtain the following non-existence result.

\begin{Prop}
Let \((M,g)\) be a closed Riemannian manifold and suppose that
the map \(q\colon M\to\s^n\subset\R^{n+1}\) defined in \eqref{eq:definition-q-biharmonic}
is a conformal-biharmonic map. If
\begin{align*}
\Scal^M>0,\qquad 
\sin^2\alpha<\frac{1}{2},\qquad \Scal^M g>3\Ric^M
\end{align*}
then the map \(q\) is constant.
\end{Prop}

\begin{proof}
Testing \eqref{eq:w-cbi} with \(v\) we obtain the following constraint
\begin{align*}
\Delta|\nabla v|^2=&(1-2\sin^2\alpha)|\nabla v|^4+\frac{2}{3}\Scal^M|\nabla v|^2
-2\Ric(\nabla v,\nabla v).
\end{align*}
Using the assumptions we then find that \(\Delta|\nabla v|^2>0\) and the maximum principle implies that \(|\nabla v|^2=const\).
Reinserting into the above equation then yields that \(|\nabla v|^2=0\)
which implies that \(v\) is constant.
\end{proof}

Note that the condition \(\Scal^M g>3\Ric^M\) implies that \(m>3\).

\subsection{The case of maps between spheres}
Since the equation for conformal-biharmonic maps involves the Ricci and the scalar curvature of the domain it will, in general, be difficult to check if maps of the form
\eqref{eq:definition-q-biharmonic} can be rendered conformal-biharmonic. Hence, we will restrict ourselves to the case of a spherical domain as this is, on the one hand, the first natural candidate to consider, and, on the other hand, already points out a
number of interesting differences to the case of standard biharmonic maps.

Thus, let us assume that \(M=\s^m\) with the round metric. Then,
we have \(\Ric^{\s^m}=(m-1)g\) and \(\Scal^{\s^m}=m(m-1)\) such that
the conformal bienergy \eqref{eq:bienergy-conformal} 
for a map \(\phi\colon\s^m\to\s^n\)
acquires the form
\begin{align*}
E_2^c(\phi)=\frac{1}{2}\int_{\s^m}\big(|\tau(\phi)|^2+\frac{2}{3}(m-1)(m-3)|d\phi|^2\big)\dv.    
\end{align*}
In the extrinsic picture, 
for \(u\colon\s^m\to\s^n\subset\R^{n+1}\), this amounts to
\begin{align*}
E_2^c(u)=\frac{1}{2}\int_M\big(|\Delta u|^2-|\nabla u|^4
+\frac{2}{3}(m-1)(m-3)|\nabla u|^2\big)\dv    
\end{align*}
and the Euler-Lagrange equation
\eqref{eq:euler-lagrange-cbi} simplifies to
\begin{align}
\label{eq:c-biharmonic-u}
\Delta^2u&+2\operatorname{div}(|\nabla u|^2\nabla u)    
-\frac{2}{3}(m-1)(m-3)\Delta u \\
\nonumber&-\big(\langle\Delta^2u,u\rangle-2|\nabla u|^4+\frac{2}{3}(m-1)(m-3)|\nabla u|^2\big)u=0.
\end{align}
Note that in the case of a spherical domain harmonic maps provide a class of solutions
to the equations for conformal-biharmonic maps and, as for biharmonic maps, it becomes important
to find the non-harmonic solutions of \eqref{eq:c-biharmonic-u}.

For the sake of completeness we want to mention that conformal-biharmonic maps between spheres are very similar to interpolating sesqui-harmonic maps to spheres which were studied in detail in \cite{MR4142862}.

In the following we will construct a number of particularly simple solutions
to the above equation.

\begin{Prop}
A map \(q\colon\s^m\to \s^n\subset\R^{n+1}\) of the form 
\eqref{eq:definition-q-biharmonic}
is conformally biharmonic if and only if    
\begin{align}
\label{eq:conformal-biharmonic-a}
(\Delta|\nabla v|^2)v=2\nabla(|\nabla v|^2)\nabla v+|\nabla v|^4v
+2\sin^2\alpha\operatorname{div}(|\nabla v|^2\nabla v)
+\frac{2}{3}(m-1)(m-3)|\nabla v|^2v.
\end{align}
\end{Prop}

\begin{proof}
This follows from a direct calculation using \eqref{eq:definition-q-biharmonic} and the fact that \(v\) is harmonic such that we can apply \eqref{eq:v-harmonic-laplace}.
\end{proof}

Testing \eqref{eq:conformal-biharmonic-a} with \(v\) we get the following constraint
\begin{align*}
\Delta|\nabla v|^2=(1-2\sin^2\alpha)|\nabla v|^4 + \frac{2}{3}(m-1)(m-3)|\nabla v|^2
\end{align*}
and in sharp contrast to the case of (standard) biharmonic maps, in general,
the maximum principle no longer forces \(|\nabla v|^2\) to be constant.

However, the above constraint immediately leads to the following restrictions
which we summarize by the following

\begin{Prop}
Let \(q\colon\s^m\to \s^n\subset\R^{n+1}\) be a conformal-biharmonic map of the form \eqref{eq:definition-q-biharmonic}. 
    \begin{enumerate}
    \item If \(m=1,3\), then \(\alpha=\frac{\pi}{4}\) and \(|\nabla v|^2=const\) and we are back to the case covered by Theorem \ref{thm:biharmonic-w}.
    \item If \(m=2\) and \(1-2\sin^2\alpha<0\), then the map \(q\) is constant.
    \item If \(m\geq 4\) and \(1-2\sin^2\alpha>0\), then the map \(q\) is constant.
\end{enumerate}
\end{Prop}

\begin{proof}
This follows from an application of the maximum principle.
\end{proof}

On the other hand, we get a large family of solutions of the equation for 
conformal-biharmonic maps between spheres by assuming that
\(v\) is an eigenmap, that is \(|\nabla v|^2=\lambda\) is constant.

\begin{Theorem}
\label{thm:existence-c-biharmonic}
The map \(q\colon\s^m\to \s^n\subset\R^{n+1}\) given by
\begin{align*}
q=(\sin\alpha\cdot v,\cos\alpha), \qquad \alpha\in(0,\frac{\pi}{2}),
\end{align*}
where \(v\colon\s^m\to \s^{n-1}\subset\R^{n}\) is a harmonic map,
is conformal-biharmonic if 
\begin{align}
\label{eq:condition-alpha-c-biharmonic}
|\nabla v|^2=\lambda,\qquad \sin^2\alpha
=\frac{1}{3}\frac{(m-1)(m-3)}{\lambda}+\frac{1}{2}.    
\end{align}
\end{Theorem}

\begin{proof}
Using that \(|\nabla v|^2=\lambda\neq 0\) in \eqref{eq:conformal-biharmonic-a}
we obtain 
\begin{align*}
\lambda v+2\sin^2\alpha\Delta v+\frac{2}{3}(m-1)(m-3) v=0    
\end{align*}
and since \(v\) is harmonic and due to our choice of \(\alpha\) we obtain the desired result.
\end{proof}

A direct calculation shows that the maps constructed in the above Theorem
have constant conformal bienergy, i.e.
\begin{align*}
E_2^c(q)=\frac{1}{2}\vol(\s^m)\big(\frac{1}{4}\lambda^2+\frac{1}{3}(m-1)(m-3)\lambda
+\frac{1}{9}(m-1)^2(m-3)^2\big)
\end{align*}
again indicating that these are unstable critical points of the conformal bienergy,
which is confirmed by the following Theorem.

\begin{Theorem}
The conformal-biharmonic maps constructed in Theorem \ref{thm:existence-c-biharmonic} are unstable critical points of the conformal bienergy \eqref{eq:bienergy-conformal}.    
\end{Theorem}

\begin{proof}
In \cite[Theorem 4.10]{BNO} it was proved that a non-harmonic conformal-biharmonic map from an Einstein manifold to the sphere with \(|\tau(\phi)|^2=const\) is always unstable.
As \(\s^m\) clearly is an Einstein manifold and as it can easily be checked that 
\(|\tau(\phi)|^2=const\) we can conclude the proof.
\end{proof}
This completes the proof of Theorem \ref{thm:c-biharmonic-main}.

\begin{Bsp}
If we assume that \(v\) is an isometric immersion, then \(\lambda=m\)
and the condition \eqref{eq:condition-alpha-c-biharmonic} becomes
\begin{align*}
\sin^2\alpha=\frac{1}{3}\frac{(m-1)(m-3)}{m}+\frac{1}{2}.    
\end{align*}
Then, we have to require
\begin{align*}
m^2-\frac{11}{2}m+3<0.    
\end{align*}
It can easily be checked that this constraint can only be satisfied 
if \(1\leq m\leq 4\) and leads to the small conformal-biharmonic hyperspheres
obtained in \cite[Theorem 3.11]{BNO}.
\end{Bsp}

However, by assuming that \(v\colon\s^m\to\s^{n-1}\subset\R^{n}\) is an eigenmap
between spheres we can easily produce additional solutions.

\begin{Prop}
There exists an \(\alpha\in(0,\frac{\pi}{2})\) such that 
the map \(q\colon\s^m\to \s^n\subset\R^{n+1}\) given by
\begin{align*}
q=(\sin\alpha\cdot v,\cos\alpha),
\end{align*}
where \(v\colon\s^m\to \s^{n-1}\subset\R^{n}\) is 
an eigenmap of degree \(k\), is rendered conformally biharmonic.
Moreover, \(\alpha\) has to satisfy 
\begin{align*}
\sin^2\alpha=\frac{1}{3}\frac{(m-1)(m-3)}{k(k+m-1)}+\frac{1}{2}
\end{align*}
leading to the constraint
\begin{align}
\label{eq:c-biharmonic-spherical-harmonic}
k>\frac{1}{6}\big(\sqrt{3}\sqrt{11m^2-38m+27}-3m+3\big).    
\end{align}
\end{Prop}

\begin{proof}
In order to obtain a non-trivial solution of \eqref{eq:condition-alpha-c-biharmonic}
we get the constraint
\begin{align*}
\frac{1}{3}(m-1)(m-3)<\frac{\lambda}{2}
\end{align*}
from \eqref{eq:condition-alpha-c-biharmonic}.
Recall that for an eigenmap of degree \(k\) we have
\(\lambda=k(k+m-1)\), see Proposition \ref{prop:homo-pol}, such that the above inequality can be rewritten as
\begin{align*}
0<k^2+k(m-1)-\frac{2}{3}m^2+\frac{8}{3}m-2.    
\end{align*}
Solving this equation for \(k\) then completes the proof.
\end{proof}

\begin{Bem}
For given \(m\) we can always find a \(k\in\N\) such that the constraint 
\eqref{eq:c-biharmonic-spherical-harmonic} is satisfied. Hence, in any dimension we can construct a non-harmonic conformal-biharmonic map out of an eigenmap of sufficiently large degree \(k\).
\end{Bem}

\subsection{Manufacturing a conformal-biharmonic out of two harmonic maps}
As in the case of biharmonic maps, we finally consider a map 
\(w\colon\s^{m}\to\s^n\subset\R^{n+1}\)
(as defined in \eqref{eq:intro-ansatz-product})
for which we make an ansatz of the form
\begin{align*}
w=(\sin\beta\cdot v_1,\cos\beta\cdot v_2),\qquad \beta\in (0,\frac{\pi}{2}),
\end{align*}
where \(v_i\colon\s^{m}\to\s^{n_i},i=1,2\) are two given harmonic maps
and \(n_1+n_2=n-1\).

The next Proposition gives a characterization when maps of the form 
\eqref{eq:intro-ansatz-product} are conformal-biharmonic.
\begin{Prop}
The map \eqref{eq:intro-ansatz-product} is conformal-biharmonic if
and only if
\begin{align}
\label{eq:equation-c-biharmonic-product}
\big(\Delta|\nabla v_i|^2\big)v_i=&-2(\nabla|\nabla v_i|^2)\nabla v_i    
+|\nabla v_i|^4v_i
+2\nabla(\sin^2\beta|\nabla v_1|^2+\cos^2\beta|\nabla v_2|^2)\nabla v_i \\
\nonumber &+(\frac{2}{3}(m-1)(m-3)-2\sin^2\beta|\nabla v_1|^2-2\cos^2\beta|\nabla v_2|^2)|\nabla v_i|^2v_i \\
\nonumber &+\big(\sin^2\beta(\Delta|\nabla v_1|^2-|\nabla v_1|^4)
+\cos^2\beta(\Delta|\nabla v_2|^2-|\nabla v_2|^4) \\
\nonumber &+2\sin^4\beta|\nabla v_1|^4+2\cos^4\beta|\nabla v_2|^4
+4\sin^2\beta\cos^2\beta|\nabla v_1|^2|\nabla v_2|^2 \\
\nonumber&-\frac{2}{3}(m-1)(m-3)(\sin^2\beta|\nabla v_1|^2+\cos^2\beta|\nabla v_2|^2)
\big)v_i,
\end{align}
where \(i=1,2\).
\end{Prop}

\begin{proof}
This follows from a direct calculation using that \(v_i,i=1,2\) are harmonic maps.
\end{proof}

\begin{Prop}
Consider a map of the form \eqref{eq:intro-ansatz-product}.
In order to be conformal-biharmonic the following constraint needs to be satisfied
\begin{align}
\label{eq:constraint-product-c-biharmonic}
\Delta\big(|\nabla v_1|^2-|\nabla v_2|^2\big)
=&\cos 2\beta(|\nabla v_1|^2-|\nabla v_2|^2)^2 \\
\nonumber&+\frac{2}{3}(m-1)(m-3)(|\nabla v_1|^2-|\nabla v_2|^2).
\end{align}    
\end{Prop}

\begin{proof}
This follows from 
testing \eqref{eq:equation-c-biharmonic-product} with \(v_i,i=1,2\) 
and the application of some trigonometric identities as in the proof of Proposition \ref{prop:product-biharmonic-constraint}.
\end{proof}

\begin{Theorem}
\label{thm:existence-c-prod-biharmonic}
Consider a map of the form \eqref{eq:intro-ansatz-product},
where \(v_i\colon\s^{m}\to\s^{n_i}\subset\R^{n_i+1},i=1,2\) are both eigenmaps
with \(|\nabla v_i|^2=\lambda_i,i=1,2\).
Then, \(w\) is a non-harmonic conformal-biharmonic map if and only if
\begin{align}
 \label{eq:product-c-biharmonic-alpha}
\lambda_1\neq\lambda_2, \qquad 
\cos 2\beta=-\frac{2}{3}(m-1)(m-3)\frac{1}{\lambda_1-\lambda_2}.
\end{align}
In particular, if \(v_i\colon\s^{m}\to\s^{n_i}\subset\R^{n_i+1},i=1,2\) are both eigenmaps
of degree \(k_i,i=1,2\) we obtain the constraint
\begin{align*}
\frac{2}{3}(m-1)(m-3)<(k_1-k_2)(k_1+k_2-m-1).
\end{align*}
\end{Theorem}

\begin{proof}
In can easily be checked that the above combination of parameters 
provides a solution to both the constraint \eqref{eq:constraint-product-c-biharmonic}
and the Euler-Lagrange equation \eqref{eq:equation-c-biharmonic-product}.

If we assume that both \(v_i\colon\s^{m}\to\s^{n_i}\subset\R^{n_i+1},i=1,2\) are eigenmaps of degree \(k_i\), then we have 
\(\lambda_i=k_i(k_i+m-1)\). In this case
\eqref{eq:product-c-biharmonic-alpha} leads to 
\begin{align*}
\cos 2\beta=
          -\frac{2}{3}\big(m-1\big)\big(m-3\big)
          \frac{1}{(k_1-k_2)(k_1+k_2-m-1)}
\end{align*}
which implies the constraint presented above.
\end{proof}

\begin{Bem}
There might also exist solutions to the constraint \eqref{eq:constraint-product-c-biharmonic} where \(|\nabla v_i|^2,i=1,2\) are non-constant but we do not further investigate this aspect here.
\end{Bem}

A direct application of \cite[Theorem 4.10]{BNO} now gives

\begin{Theorem}
The conformal-biharmonic maps constructed in Theorem \ref{thm:existence-c-prod-biharmonic} are unstable critical points of the conformal bienergy \eqref{eq:bienergy-conformal}.    
\end{Theorem}

This completes the proof of Theorem \ref{thm:c-biharmonic-product-main}.

\bibliographystyle{plain}
\bibliography{mybib}
\end{document}